figure 2, curves $\chi^n$ which reach $(n,0)$ and afterwards are horizontal until $(n+1/4,0)$ where they end into $T_n$. For every $n$ there is a curve $\psi^n$ which starts from $(n+3/4,0)$ and has $t = n+1$ as vertical asymptote to $-\infty$. For every $n$ and $i$ there is a curve $\varphi_i^n$ which starts from $(n+3/4, 2^{-i+1})$ on the right side of $T_n$ and reaches $(n+1,0)$, where it merges with $\chi^{n+1}$.

Obviously all these curves should be arranged so that at their intersections they have the same slope and also so that the function describing their derivatives is continuous: this in particular can be obtained by requiring that $\varphi_i^n$ keeps close to $\psi^n$ for longer and longer portions of the latter as $i$ gets larger.

Now we define $F_\alpha$ on $\mathbb{R}^2 \setminus T$ by requiring that all these curves are (portions of) solutions of the differential equation given by $F_\alpha$. Again we would like to avoid the existence of solutions of that differential equation which passes through one point on the curve $\psi^n$ and then moves to the right of it. Similarly there should be no points on the $\chi^n$'s or the $\varphi_i^n$'s where local multiplicity occurs towards the right. We furthermore want that $F_\alpha(t,0) = 0$ and $F_\alpha(t,1) = 1$ for every $n$ and $t \in [n+1/4, n+3/4]$ and that $F_\alpha(n+j/4, t) \geq 0$ for every $n$, $j = 1, 3$ and $t \in [0,1]$. It should also exist a solution that moves from $(-1, 0)$ to the left to $-\infty$. This can be attained more or less in the same way we attained the corresponding requirements in the proof of theorem 5.6.

The definition of $F_\alpha$ within each $T_n$ follows again the same ideas employed in the proof of theorem 5.6 and we are left with showing that this construction actually works.

First suppose that $\alpha$ is such that $\forall^\infty n \; \forall^\infty m \; \alpha(n,m) = 0$. This means that for all but finitely many $n$'s we have local multiplicity within $T_n$. Let $n_0$ be such that for every $n \geq n_0$ this happens. From $(-1,0)$ follow $\chi$ until $\chi^{n_0}$ splits and follow the latter to $(n_0 + 1/4, 0)$. Now avoid the possibility of being led to an asymptote by $\psi^{n_0}$ by exiting from $T_{n_0}$ above some $\varphi_i^{n_0}$. This can be repeated for all the following $T_n$'s thereby constructing a global solution of the Cauchy problem.

If $\alpha$ is such that for infinitely many $n$'s we have local uniqueness within $T_n$, no matter which $\chi^{n_0}$ we will pick to avoid being led to a vertical asymptote by $\chi$ we will be led to some other vertical asymptote by $\psi^n$ where $n \geq n_0$ is the first such that uniqueness holds in $T_n$. □

Notice that lemma 5.7 and theorem 5.6 show that the complexities of $\mathcal{G}$ and $\mathcal{G}_\forall$ are different. Contrast this with the fact that the complexities of $\mathcal{U}$ and $\mathcal{U}_\forall$ are the same (theorems 3.12 and 3.15).

Our proof shows also that for every $(x,y) \in \mathbb{R}^2$ $\mathcal{G}_{(x,y)}$ is $\mathbf{\Sigma}_4^0$-hard.

The obstacle to the existence of a global solution of the Cauchy problem whose solutions are depicted in figure 2 is that $\lim_{n \to \infty} \varphi_n(0) = -\infty$, i.e. the $\varphi_n$'s drop arbitrarily low. The next lemma shows that if this does not happen then a global solution does exist. Let us first indroduce some helpful notation and prove an easy proposition about it.

**Definition 5.8.** If $F \in C(\mathbb{R}^2)$ and $(x,y) \in \mathbb{R}^2$ let
$$\mathcal{S}_{F,(x,y)}^+ = \left\{ \varphi \upharpoonright [x, +\infty) \mid \varphi \in \mathcal{S}_{F,(x,y)} \right\}.$$

In other words $\mathcal{S}_{F,(x,y)}^+$ is the set of $\varphi \in C^1([x,b))$ for some $b$ such that $x < b \leq +\infty$ satisfying $\varphi(x) = y$, $\varphi'(t) = F(t, \varphi(t))$ for every $t \in [x,b)$ (in $x$ we are considering the right derivative), and which are non-extendible to the right.

**Definition 5.9.** If $\{\varphi_n\}_n$ is a sequence of real-valued functions defined on subsets of $\mathbb{R}$ let $\inf_n \varphi_n$ be the function defined on a subset of $\bigcup_n \text{dom}(\varphi_n)$ by letting $\inf_n \varphi_n(t) = \inf \{ \varphi_n(t) \mid t \in \text{dom}(\varphi_n) \}$ (the function is defined only if the inf is greater then $-\infty$).

**Definition 5.10.** If $\varphi$ is a real-valued function defined on an interval in $\mathbb{R}$ and $b = \sup(\mathrm{dom}(\varphi))$ we denote $\lim_{t \to b^-} \varphi(t)$ by $\lim_+ \varphi$.

**Definition 5.11.** If $\varphi$ and $\psi$ are two real-valued functions defined on subsets of $\mathbb{R}$ we write $\varphi \leq \psi$ if $\varphi(t) \leq \psi(t)$ whenever $t \in \mathrm{dom}(\varphi) \cap \mathrm{dom}(\psi)$. We also define two functions $\varphi \vee \psi$ and $\varphi \wedge \psi$ with domain $\mathrm{dom}(\varphi) \cup \mathrm{dom}(\psi)$. Both functions coincide with $\varphi$ on $\mathrm{dom}(\varphi) \setminus \mathrm{dom}(\psi)$ and with $\psi$ on $\mathrm{dom}(\psi) \setminus \mathrm{dom}(\varphi)$. If $t \in \mathrm{dom}(\varphi) \cap \mathrm{dom}(\psi)$ then $(\varphi \vee \psi)(t) = \max(\varphi(t), \psi(t))$ and $(\varphi \wedge \psi)(t) = \min(\varphi(t), \psi(t))$.

**Proposition 5.12.** *If $F \in C(\mathbb{R}^2)$, $(x, y) \in \mathbb{R}^2$ and $\varphi, \psi \in \mathcal{S}^+_{F,(x,y)}$ are such that $\mathrm{dom}(\varphi) \subseteq \mathrm{dom}(\psi)$ then:*

(1) $\lim_+ \varphi \neq -\infty$ *implies* $\varphi \wedge \psi \in \mathcal{S}^+_{F,(x,y)}$;

(2) $\lim_+ \varphi \neq +\infty$ *implies* $\varphi \vee \psi \in \mathcal{S}^+_{F,(x,y)}$.

*Proof.* The condition on $\lim_+ \varphi$ insures that $\varphi \wedge \psi$ (resp. $\varphi \vee \psi$) is continuous. The fact that whenever $\varphi(t) = \psi(t)$ we have also $\varphi'(t) = \psi'(t)$ immediately yields the conclusion. □

The next lemma gives a way to build new solutions of an ODE out of old ones.

**Lemma 5.13.** *Let $(F, x, y) \in C(\mathbb{R}^2) \times \mathbb{R}^2$ and let $\{\varphi_n\}_n$ be a sequence of elements of $\mathcal{S}^+_{F,(x,y)}$ such that $\lim_+ \varphi_n \neq -\infty$ for every $n$. Let $T \in (x, +\infty]$ be such for every $t \in [x, T) \subseteq \bigcup_n \mathrm{dom}(\varphi_n)$ we have $\inf\{\varphi_n(t) \mid t \in \mathrm{dom}(\varphi_n)\} > -\infty$ while, if $T < +\infty$, $\inf\{\varphi_n(T) \mid T \in \mathrm{dom}(\varphi_n)\} = -\infty$. Then $\inf_n \varphi_n \upharpoonright [x, T) \in \mathcal{S}^+_{F,(x,y)}$. In particular if $T = +\infty$ this means that $\inf_n \varphi_n \in \mathcal{S}^+_{F,(x,y)}$.*

*Proof.* Replacing $\varphi_n$ by $\varphi_0 \wedge \ldots \wedge \varphi_n$ (and using proposition 5.12) we may assume that $\mathrm{dom}(\varphi_n) \subseteq \mathrm{dom}(\varphi_{n+1})$ and $\varphi_n \geq \varphi_{n+1}$. Therefore the graph of $\inf_n \varphi_n$ is the limit within $\mathbf{F}(\mathbb{R}^2)$ of the graphs of the $\varphi_n$'s. The conclusion follows by lemma 2.7. □

We now define two subsets of $\mathcal{G}^+$.

**Definition 5.14.** Let $\mathcal{H}_-$ be the set of $(F, x, y) \in C(\mathbb{R}^2) \times \mathbb{R}^2$ such that there exist $N, h \in \mathbb{N}$ such that for every $n \geq 1$ there exists $\varphi_n \in \mathcal{S}^+_{F,(x,y)}$ such that $[x, x+N+n] \subset \mathrm{dom}(\varphi_n)$, $\|\varphi_n \upharpoonright [x, x+N]\|_\infty \leq h$, and there is no $\varphi \in \mathcal{S}^+_{F,(x,y)}$ such that $[x, x+N] \subset \mathrm{dom}(\varphi)$, $\|\varphi \upharpoonright [x, x+N]\|_\infty \leq h$, and $\lim_+ \varphi = -\infty$. In other words $(F, x, y) \in \mathcal{H}_-$ if and only if there exists a rectangle of base $[x, x+N]$ and height $[-h, h]$ such that the Cauchy problem admits partial solutions of arbitrary length which are trapped inside the rectangle $[x, x+N] \times [-h, h]$, but no solutions trapped inside the rectangle can have a vertical asymptote going to $-\infty$.

$\mathcal{H}_+$ is defined analogously asking that no solutions satisfying the rectangle condition go to $+\infty$ at their right vertical asymptote

**Lemma 5.15.** $\mathcal{H}_- \cup \mathcal{H}_+ \subseteq \mathcal{G}^+$.

*Proof.* The symmetry of the definitions of $\mathcal{H}_-$ and $\mathcal{H}_+$ entails that it suffices to prove $\mathcal{H}_- \subseteq \mathcal{G}^+$. Let $(F, x, y) \in \mathcal{H}_-$ and let $N, h$ witness this. Moreover let $\{\varphi_n\}_n$ be a sequence of witnesses for the first part of the definition.

We claim that $\inf\{\varphi_n(t) \mid t \in \mathrm{dom}(\varphi_n)\} > -\infty$ for every $t \in [x, +\infty)$. If $t \in [x, x+N]$ then this follows from the fact that $\varphi_n(t) \geq -h$ for every $n$. Suppose now, towards a contradiction, that $t > x+N$ is such that $\inf\{\varphi_n(t) \mid t \in \mathrm{dom}(\varphi_n)\} = -\infty$. Then lemma 5.13 implies that for some $T \leq t$, $\varphi = \inf_n \varphi_n \upharpoonright [x, T)$ is an element of $\mathcal{S}^+_{F,(x,y)}$ satisfying $[x, x+N] \subset \mathrm{dom}(\varphi)$, $\|\varphi \upharpoonright [x, x+N]\|_\infty \leq h$, and $\lim_+ \varphi = -\infty$. This contradicts $(F, x, y) \in \mathcal{H}_-$ and proves the claim.

The claim and lemma 5.13 imply that $(F, x, y) \in \mathcal{G}^+$. □

We now define another subset of $\mathcal{G}^+$.

**Definition 5.16.** Let $\mathcal{H}_\infty$ be the set of $(F, x, y) \in C(\mathbb{R}^2) \times \mathbb{R}^2$ such that for every $n \geq 1$ there exist $\varphi_n, \psi_n \in \mathcal{S}^+_{F,(x,y)}$ such that $[x, x+n] \subset \text{dom}(\varphi_n) \cap \text{dom}(\psi_n)$, $\lim_+ \varphi_n = +\infty$, and $\lim_+ \psi_n = -\infty$. In other words $(F, x, y) \in \mathcal{H}_\infty$ if and only if the Cauchy problem admits both partial solutions of arbitrary length such that $\lim_+ \varphi = -\infty$ and partial solutions of arbitrary length such that $\lim_+ \varphi = +\infty$.

**Lemma 5.17.** $\mathcal{H}_\infty \subseteq \mathcal{G}^+$.

*Proof.* Suppose $(F, x, y) \in \mathcal{H}_\infty$ and let $\{\varphi_n\}_n$ and $\{\psi_n\}_n$ be sequences that witness this. If the domain of some $\varphi_n$ or $\psi_n$ is unbounded above we are done. Otherwise let $b_n = \sup(\text{dom}(\varphi_n))$ and $d_n = \sup(\text{dom}(\psi_n))$. We can suppose that $b_n \leq d_n < b_{n+1}$ for every $n$. We define inductively two sequences $\{\tilde\varphi_n\}_n$ and $\{\tilde\psi_n\}_n$ of elements of $\mathcal{S}^+_{F,(x,y)}$ as follows:

$$\tilde\varphi_0 = \varphi_0; \qquad\qquad \tilde\psi_0 = \varphi_0 \wedge \psi_0;$$
$$\tilde\varphi_{n+1} = \tilde\psi_n \vee (\tilde\varphi_n \wedge \varphi_{n+1}); \qquad \tilde\psi_{n+1} = \tilde\varphi_{n+1} \wedge (\tilde\psi_n \vee \psi_{n+1}).$$

A straightforward induction shows that for every $n$ $\tilde\varphi_n, \tilde\psi_n \in \mathcal{S}^+_{F,(x,y)}$ (using proposition 5.12), $\tilde\varphi_n \geq \tilde\varphi_{n+1} \geq \tilde\psi_{n+1} \geq \tilde\psi_n$, $\text{dom}(\tilde\varphi_n) = \text{dom}(\varphi_n)$, $\text{dom}(\tilde\psi_n) = \text{dom}(\psi_n)$, $\lim_+ \tilde\varphi_n = +\infty$, and $\lim_+ \tilde\psi_n = -\infty$.

If $t > x$ let $N$ be such that $t < x + N$: for every $n$ such that $t \in \text{dom}(\tilde\varphi_n)$ we have $\tilde\varphi_n(t) \geq \tilde\psi_N(t)$. Hence $\inf\{\tilde\varphi_n(t) \mid t \in \text{dom}(\tilde\varphi_n)\} > -\infty$ and we can apply lemma 5.13 to obtain that $\tilde\varphi = \inf_n \tilde\varphi_n \in \mathcal{S}^+_{F,(x,y)}$. As $\text{dom}(\tilde\varphi) = \bigcup_n \text{dom}(\tilde\varphi_n) \supseteq \bigcup_n \text{dom}(\varphi_n) \supseteq \bigcup_n [x, x+n]$ it follows that $\text{dom}(\tilde\varphi)$ is unbounded above. Therefore $(F, x, y) \in \mathcal{G}^+$. $\square$

$\mathcal{H}_- \cup \mathcal{H}_+$ and $\mathcal{H}_\infty$ show two different ways in which a Cauchy problem can have a solution with domain unbounded above. These are actually the only possible ways this can happen.

**Theorem 5.18.** $\mathcal{G}^+ = \mathcal{H}_- \cup \mathcal{H}_+ \cup \mathcal{H}_\infty$.

*Proof.* Lemmas 5.15 and 5.17 show that $\mathcal{H}_- \cup \mathcal{H}_+ \cup \mathcal{H}_\infty \subseteq \mathcal{G}^+$, so we need to prove only the reverse inclusion.

Let $(F, x, y) \in \mathcal{G}^+$ and define:

$$B^+ = \left\{ b \mid \exists \varphi \in \mathcal{S}^+_{F,(x,y)} (\text{dom}(\varphi) = [x, b] \,\&\, \lim_+ \varphi = +\infty) \right\} \cup \{x\}$$
$$B^- = \left\{ b \mid \exists \varphi \in \mathcal{S}^+_{F,(x,y)} (\text{dom}(\varphi) = [x, b] \,\&\, \lim_+ \varphi = -\infty) \right\} \cup \{x\}$$

and set

$$b^+ = \sup B^+ \qquad \text{and} \qquad b^- = \sup B^-.$$

If $b^+ = b^- = +\infty$ then clearly $(F, x, y) \in \mathcal{H}_\infty$. Now we will show that if $b^- < +\infty$ then $(F, x, y) \in \mathcal{H}_-$: a symmetric argument shows that if $b^+ < +\infty$ then $(F, x, y) \in \mathcal{H}_+$ and completes the proof of the theorem.

Suppose $b^- < +\infty$ and let $N \in \mathbb{N}$ be such that $b^- < x + N$. Since $(F, x, y) \in \mathcal{G}^+$ there exists $\psi \in \mathcal{S}^+_{F,(x,y)}$ such that $\text{dom}(\psi) = [x, +\infty)$. Fix such a $\psi$ and let $h \in \mathbb{N}$ be such that $\|\psi \restriction [x, x+N]\|_\infty \leq h$. We claim that $N$ and $h$ witness $(F, x, y) \in \mathcal{H}_-$. For every $n$ we can set $\varphi_n = \psi$ so that the first part of the definition (the existence of solutions of arbitrary length) is satisfied. On the other hand the choice of $N$ guarantees that also the second part of the definition (the non existence of solutions with $\lim_+ \varphi = -\infty$ which are initially contained in the rectangle $[x, x+N] \times [-h, h]$) is satisfied. $\square$

Theorem 5.18 shows that to find an upper bound for the complexity of $\mathcal{G}^+$ (and hence of $\mathcal{G}$) it suffices to find upper bounds for the complexities of $\mathcal{H}_-$ and $\mathcal{H}_\infty$ (obviously the complexity of $\mathcal{H}_+$ is the same of that of $\mathcal{H}_-$). This is precisely what we are going to do, using the following technical lemma.

**Lemma 5.19.** *(a) Fix $N, h \in \mathbb{N}$. Let $A_{N,h}^+$ be the set of all Cauchy problems which have a solution with a vertical asymptote to $+\infty$ at some $\xi > x + N$ and such that this solution is bounded by $h$ on the interval $[x, x+N]$, i.e.*

$$A_{N,h}^+ = \{\, (F, x, y) \in C(\mathbb{R}^2) \times \mathbb{R}^2 \mid \exists \varphi \in \mathcal{S}_{F,(x,y)}^+ [\exists \xi > x + N (\operatorname{dom}(\varphi) = [x, \xi)) \\ \&\, \lim\nolimits_+ \varphi = +\infty\, \&\, \|\varphi \restriction [x, x+N]\|_\infty \leq h]\}.$$

*Similarly we define $A_{N,h}^-$ by replacing "$\lim_+ \varphi = +\infty$" with "$\lim_+ \varphi = -\infty$". Then both $A_{N,h}^+$ and $A_{N,h}^-$ are $\mathbf{\Sigma}_2^0$.*

*(b) Fix $N \in \mathbb{N}$. Let $A_N^+$ be the set of all Cauchy problems which have a solution with a vertical asymptote to $+\infty$ at some $\xi > x + N$, i.e.*

$$A_N^+ = \{\, (F, x, y) \in C(\mathbb{R}^2) \times \mathbb{R}^2 \mid \exists \varphi \in \mathcal{S}_{F,(x,y)}^+ [\exists \xi > x + N (\operatorname{dom}(\varphi) = [x, \xi)) \\ \&\, \lim\nolimits_+ \varphi = +\infty]\}.$$

*Similarly we define $A_N^-$ by replacing "$\lim_+ \varphi = +\infty$" with "$\lim_+ \varphi = -\infty$". Then both $A_N^+$ and $A_N^-$ are $\mathbf{\Sigma}_2^0$.*

*(c) Fix $M, N, h \in \mathbb{N}$, with $M \geq N$. Let $B_{N,h}^M$ be the set of all Cauchy problems which have solutions of length at least $M$ which are bounded by $h$ on the interval $[x, x+N]$, i.e.*

$$B_{N,h}^M = \{\, (F, x, y) \in C(\mathbb{R}^2) \times \mathbb{R}^2 \mid \exists \varphi \in \mathcal{S}_{F,(x,y)}^+ [\operatorname{dom}(\varphi) \supseteq [x, x+M] \\ \&\, \|\varphi \restriction [x, x+N]\|_\infty \leq h]\}.$$

*Then $B_{N,h}^M$ is $\mathbf{\Sigma}_2^0$.*

*Proof.* (a) We claim that

(*) $$A_{N,h}^+ = \bigcup_{b,\ell \in \mathbb{Q}} \bigcap_{M \in \mathbb{N}} D_{b,\ell}^M$$

where

$$D_{b,\ell}^M = \{\, (F, x, y) \in C(\mathbb{R}^2) \times \mathbb{R}^2 \mid \exists r \in \mathbb{R}\, \exists \varphi \in \mathcal{S}_{F,(x,y)}^+ [x + N \leq r \leq x + N + \ell \\ \&\, \|\varphi \restriction [x, x+N]\|_\infty \leq h\, \&\, r \in \operatorname{dom}(\varphi)\, \&\, (\forall t \in [x, r]\, \varphi(t) \geq b)\, \&\, \varphi(r) \geq M]\}.$$

In fact if $\varphi$ witnesses $(F, x, y) \in A_{N,h}^+$ let $\operatorname{dom}(\varphi) = [x, \xi)$ with $\xi > x + N$ and recall that $\lim_+ \varphi = +\infty$. Then we can pick $b, \ell \in \mathbb{Q}$ such that $\ell > 0$, $x + N + \ell \geq \xi$ and $b < \inf\{\varphi(t) \mid x \leq t \leq \xi\}$. Therefore for every $M \in \mathbb{N}$ there exists $r$ that, together with $\varphi$, witnesses $(F, x, y) \in D_{b,\ell}^M$.

Conversely suppose $(F, x, y) \in \bigcap_{M \in \mathbb{N}} D_{b,\ell}^M$ for some $b, \ell \in \mathbb{Q}$. Then for every $M \in \mathbb{N}$ there exist $r_M$ and $\varphi_M$ satisfying the conditions in the definition of $D_{b,\ell}^M$. As $[x+N, x+N+\ell]$ and $\mathbf{F}(\mathbb{R}^2)$ (with the Fell topology) are compact, by extracting a subsequence we may assume that $\lim_{M \to +\infty} r_M = r \in [x+N, x+N+\ell]$ and $\lim_{M \to \infty} \varphi_M = \Gamma \in \mathbf{F}(\mathbb{R}^2)$. By lemma 2.7 there exists $\varphi \in \mathcal{S}_{F,(x,y)}^+$ whose graph is $(\operatorname{dom}(\varphi) \times \mathbb{R}) \cap \Gamma$ and such that $[x, x+N] \subset \operatorname{dom}(\varphi) \subseteq [x, r)$, $\|\varphi \restriction [x, x+N]\|_\infty \leq h$ and $\varphi$ is bounded below by $b$. Therefore $\varphi$ has a vertical asymptote to $+\infty$ at some $\xi \leq r$, and hence $(F, x, y) \in A_{N,h}^+$.

Notice that in proving (*) we actually verified that $D_{b,\ell}^M$ is closed and therefore $A_{N,h}^+$ is $\mathbf{\Sigma}_2^0$. The result about $A_{N,h}^-$ is completely analogous.

(b) follows from (a) and the equalities $A_N^+ = \bigcup_{h \in \mathbb{N}} A_{N,h}^+$ and $A_N^- = \bigcup_{h \in \mathbb{N}} A_{N,h}^-$.

(c) The argument is similar to the one used to prove (a), only simpler. Fix $M$, $N$ and $h$. Then
$$B_{N,h}^M = \bigcup_{\ell \in \mathbb{Q}^+} E_\ell$$
with
$$E_\ell = \{\, (F, x, y) \in C(\mathbb{R}^2) \times \mathbb{R}^2 \mid \exists \varphi \in \mathcal{S}_{F,(x,y)}^+ [\mathrm{dom}(\varphi) \supset [x, x+M]$$
$$\&\, \|\varphi \restriction [x, x+M]\|_\infty \le \ell \,\&\, \|\varphi \restriction [x, x+N]\|_\infty \le h]\,\}.$$

By Ascoli-Arzelà or, equivalently, by the compactness of $\mathbf{F}(\mathbb{R}^2)$, each $E_\ell$ is closed and hence $B_{N,h}^M$ is $\boldsymbol{\Sigma}_2^0$. □

**Lemma 5.20.** $\mathcal{H}_-$ *and* $\mathcal{H}_+$ *are* $\boldsymbol{\Sigma}_4^0$. $\mathcal{H}_\infty$ *is* $\boldsymbol{\Pi}_3^0$.

*Proof.* Using the notations of lemma 5.19 we have
$$\mathcal{H}_- = \bigcup_{h, N \in \mathbb{N}} \left[ \left( (C(\mathbb{R}^2) \times \mathbb{R}^2) \setminus A_{N,h}^- \right) \cap \bigcap_{M \ge N} B_{N,h}^M \right],$$
$$\mathcal{H}_+ = \bigcup_{h, N \in \mathbb{N}} \left[ \left( (C(\mathbb{R}^2) \times \mathbb{R}^2) \setminus A_{N,h}^+ \right) \cap \bigcap_{M \ge N} B_{N,h}^M \right],$$
and
$$\mathcal{H}_\infty = \bigcap_{M \in \mathbb{N}} \left( A_M^+ \cap A_M^- \right).$$

Therefore $\mathcal{H}_-$ and $\mathcal{H}_+$ are $\boldsymbol{\Sigma}_4^0$, while $\mathcal{H}_\infty$ is $\boldsymbol{\Pi}_3^0$. □

By modifying the proof of lemma 5.7 (notice that in that proof if $\alpha \in S_4$ then $(F_\alpha, -1, 0) \in \mathcal{H}_+$) it is not difficult to show that $\mathcal{H}_\infty$ is actually $\boldsymbol{\Pi}_3^0$-hard, and hence $\boldsymbol{\Pi}_3^0$-complete. It suffices to delete $\chi$ and the various $\chi^n$'s with $n > 0$ and flip the behaviour of $F_\alpha$ inside the $T_n$'s with $n$ odd so that $\lim_+ \psi^{2k+1} = +\infty$. This construction yields $(F_\alpha, -1, 0) \in \mathcal{H}_\infty$ if and only if $\forall n\, \forall^\infty m\, \alpha(n, m) = 0$.

**Theorem 5.21.** $\mathcal{G}$ *is* $\boldsymbol{\Sigma}_4^0$-*complete.*

*Proof.* The $\boldsymbol{\Sigma}_4^0$-hardness of $\mathcal{G}$ was proved in lemma 5.7. Since the finite union of $\boldsymbol{\Sigma}_4^0$ and $\boldsymbol{\Pi}_3^0$ sets is $\boldsymbol{\Sigma}_4^0$ theorem 5.18 and lemma 5.20 show that $\mathcal{G}^+$ is $\boldsymbol{\Sigma}_4^0$. But then $\mathcal{G}^-$ is also $\boldsymbol{\Sigma}_4^0$ and $\mathcal{G} = \mathcal{G}^+ \cap \mathcal{G}^-$ is $\boldsymbol{\Sigma}_4^0$. □

DIP. DI MATEMATICA, UNIVERSITÀ DI TORINO, VIA CARLO ALBERTO 10, 10123 TORINO, ITALY
*E-mail address*: {andretta,marcone}@dm.unito.it